\documentclass[11pt,twoside,reqno]{amsart}
\usepackage{times}
\usepackage{amsmath}
\usepackage{amssymb}
\usepackage{amsfonts}
\usepackage{bbm}
\usepackage{graphicx}
\usepackage{a4wide}
\usepackage{url}

\newcommand{\R}{\mathbbm{R}}

\newcommand{\RNonNeg}{\mathbbm{R}_+}

\newcommand{\Q}{\mathbbm{Q}}

\newcommand{\Z}{\mathbbm{Z}}
\newcommand{\N}{\mathbbm{N}}
\newcommand{\ints}[1]{[{#1}]}

\DeclareMathOperator{\NP}{NP}
\DeclareMathOperator{\coNP}{coNP}

\DeclareMathOperator{\transposeOp}{t}
\newcommand{\transpose}[1]{{#1}^{\transposeOp}}

\DeclareMathOperator{\identOp}{Id}
\newcommand{\ident}[1]{\identOp_{#1}}

\DeclareMathOperator{\kernelOp}{ker}
\newcommand{\kernel}[1]{\kernelOp({#1})}

\newcommand{\unitVec}[1]{\mathbbm{e}_{#1}}

\newcommand{\zeroVec}[1]{\mathbb{O}_{#1}}

\newcommand{\scalProd}[2]{\langle{#1},{#2}\rangle}

\newcommand{\row}[2]{{#1}_{{#2},\star}}
\newcommand{\col}[2]{{#1}_{\star,{#2}}}

\newcommand{\sd}[1]{\delta({#1})}
\newcommand{\qsd}[1]{\Delta({#1})}

\DeclareMathOperator{\affOp}{aff}
\newcommand{\aff}[1]{\affOp({#1})}

\DeclareMathOperator{\convOp}{conv}
\newcommand{\conv}[1]{\convOp({#1})}
\DeclareMathOperator{\coneOp}{cone}

\DeclareMathOperator{\cconeOp}{ccone}
\newcommand{\ccone}[1]{\cconeOp({#1})}

\DeclareMathOperator{\monoOp}{mono}
\newcommand{\mono}[1]{\monoOp({#1})}

\DeclareMathOperator{\hOp}{H}
\newcommand{\hyperEbene}[2]{\hOp^=({#1},{#2})}
\newcommand{\halbRaum}[2]{\hOp^{\le}({#1},{#2})}

\DeclareMathOperator{\polyOp}{P}
\newcommand{\polyLe}[2]{\polyOp^{\le}({#1},{#2})}

\DeclareMathOperator{\charConeOp}{char}
\newcommand{\charCone}[1]{\charConeOp({#1})}

\DeclareMathOperator{\linealOp}{lineal}
\newcommand{\lineal}[1]{\linealOp({#1})}
\DeclareMathOperator{\eqSetOp}{Eq}

\newcommand{\eqSetMit}[3]{\eqSetOp_{{#1}x\le{#2}}({#3})}
\DeclareMathOperator{\intHullOp}{I}
\newcommand{\intHull}[1]{{#1}_{\intHullOp}}
\DeclareMathOperator{\homogOp}{homog}
\newcommand{\homog}[1]{\homogOp({#1})}
\newcommand{\facLat}[1]{\mathcal{L}({#1})}

\newcommand{\polar}[1]{{#1}^{\circ}}
\newcommand{\ppolar}[1]{{#1}^{\circ\circ}}

\DeclareMathOperator{\inzOp}{inc}
\newcommand{\inzDir}[1]{\inzOp({#1})}
\newcommand{\inz}[1]{\inzOp({#1})}

\newcommand{\setDef}[2]{\{{#1}:{#2}\}}

\newcommand{\charVec}[1]{\chi({#1})}
\newcommand{\enc}[1]{\langle{#1}\rangle}
\newcommand{\encmax}[1]{\langle{#1}\rangle_{\max}}
\DeclareMathOperator{\outerOp}{outer}
\newcommand{\encouter}[1]{\langle{#1}\rangle^{\outerOp}_{\max}}
\DeclareMathOperator{\innerOp}{inner}
\newcommand{\encinner}[1]{\langle{#1}\rangle^{\innerOp}_{\max}}

\newcommand{\NICHT}[1]{}

\newtheorem{theorem}{Theorem}

\begin{document}

\title{Basic Polyhedral Theory}%
\author{Volker Kaibel}%
\address{Volker Kaibel}
\email{kaibel@ovgu.de}%
\thanks{Article prepared for \emph{Wiley Encyclopedia of Operations Research and Management Science}}

\date{\today}%


\maketitle

\noindent
A \emph{polyhedron} is the 
intersection 
of finitely many affine halfspaces, where an \emph{affine halfspace} is a set
\begin{equation*}
	\halbRaum{a}{\beta}=\setDef{x\in\R^n}{\scalProd{a}{x}\le\beta}
\end{equation*}
for some $a\in\R^n$ and~ $\beta\in\R$ (here,  $\scalProd{a}{x}=\sum_{j=1}^na_jx_j$ denotes the standard scalar product on~$\R^n$). Thus, every polyhedron is the 
set
\begin{equation*}
	\polyLe{A}{b}=\setDef{x\in\R^n}{Ax\le b}
\end{equation*}  
of feasible solutions to a system $Ax\le b$ of linear inequalities for some matrix $A\in\R^{m\times n}$ and some vector $b\in\R^m$. 
Clearly, all sets $\setDef{x\in\R^n}{Ax\le b, A'x=b'}$ are polyhedra as well, as the system $A'x=b'$ of linear equations is equivalent the system $A'x\le b', -A'x\le-b'$ of linear inequalities. A bounded polyhedron is called a \emph{polytope} (where \emph{bounded} means that there is a  bound which no  coordinate of any point in the polyhedron exceeds  in absolute value). 

Polyhedra are of great importance for Operations Research, because they are not only the sets of feasible solutions to Linear Programs (LP), for which we have beautiful duality results and both practically and theoretically efficient algorithms, but even the solution of (Mixed) Integer Linear Programming (MILP) problems can be reduced to linear optimization problems over polyhedra. This relationship to a large extent forms the backbone of the extremely successful story of (Mixed) Integer Linear Programming and Combinatorial Optimization over the last few decades. 

In Section~\ref{sec:geom}, we review those parts of the general theory of polyhedra that are most important with respect to optimization questions, while in Section~\ref{sec:int} we treat concepts that are particularly relevant for  Integer Programming. Most of the ``basic polyhedral theory'' today is standard textbook knowledge. In Section~\ref{sec:ref}, for some of the results we provide references to the original papers. There, we also give pointers to proofs of the theorems mentioned in the first two sections, where we mainly refer to the beautiful book by Schrijver~\cite{Sch86}. There are, of course, many other excellent  treatments of the theory of polyhedra with respect to optimization questions, e.g., in the recent survey by  
Conforti, Cornu\'{e}jols, and Zambelli~\cite{CCZ09}, in the handbook articles by Schrijver~\cite{Sch95} and  Burkard~\cite{Bur93}, as well as 
in the books by Nemhauser and Wolsey~\cite{NW88}, Gr\"otschel, Lov\`{a}sz, and Schrijver~\cite{GLS88},
Bertsimas and Weismantel~\cite{BW05}, Cook, Cuningham, Pulleyblank, and Schrijver~\cite{CCPS98}, Wolsey~\cite{Wol98}, Korte and Vygen~\cite{KV08}, or Barvinok~\cite{Bar02}. 
The books by Ziegler~\cite{Zie95} and Gr\"unbaum~\cite{Gru03} are most important sources for the general geometric and combinatorial theory of polyhedra, in particular of polytopes. We also refer to the handbook article by Gritzmann and Klee~\cite{GK93} as well as the one by Bayer and Lee~\cite{BL93}.

\section{The Geometry of Polyhedra}
\label{sec:geom}

\subsection{Some Notation}
We define $\ints{p}=\{1,\dots,p\}$ and denote by $\RNonNeg=\setDef{\alpha\in\R}{\alpha\ge 0}$ the set of nonnegative real numbers.
A \emph{submatrix} of $M\in\R^{m\times n}$ is a matrix $M_{I,J}\in\R^{I\times J}$ for some $\varnothing\ne I\subseteq\ints{m}$ and $\varnothing\ne J\subseteq\ints{n}$ formed by the rows and columns of~$M$ indexed by the elements of~$I$ and~$J$, respectively. In particular, $M_{i,j}\in\R$ is the entry in row~$i$ and column~$j$. We write $\row{M}{I}=M_{I,\ints{n}}$ and $\col{M}{J}=M_{\ints{m},J}$, in particular, $\row{M}{i}\in\R^n$ and~$\col{M}{j}\in\R^m$ are the $i$-th row and the $j$-th column of~$M$, respectively. The \emph{kernel} of $M\in\R^{m\times n}$ is $\kernel{(M)}=\setDef{x\in\R^n}{Mx=\zeroVec{}}$.
The \emph{identity matrix} $\ident{n}\in\R^{n\times n}$ has one-entries on its main diagonal and zeroes elsewhere.

For $x\in\R^n$ and $J\subseteq\ints{n}$, the vector formed by the components of~$x$ indexed by elements of~$J$ is denoted by~$x_J\in\R^J$. 
We denote by $\scalProd{x}{y}=\sum_{j=1}^nx_jy_j$ the standard scalar product of $x,y\in\R^n$. We consider, in the context of matrix multiplication, all vectors as column vectors, and use~$\transpose{(\dots)}$ to refer to the transposed matrix or vector. We denote by $\zeroVec{}\in\R^n$ the zero vector in~$\R^n$ and by~$\unitVec{i}\in\R^n$  the standard unit vector having its $i$-th component equal to one,  all other components being zero. 

Denoting by~$\sd{M}$ the set of all determinants of submatrices (formed by  arbitrary subsets of rows and columns of equal cardinality, including the empty submatrix, whose determinant is considered to be one) of a matrix $M\in\R^{m\times n}$, we define
\begin{equation*}
	\qsd{M}=\setDef{\tfrac{p}{q}}{p,q\in \sd{M}\cup (-\sd{M}), q\ne 0}\,,
\end{equation*}
and, for every finite set $\varnothing\ne V\subseteq\R^n$, we set $\qsd{V}=\qsd{M}$, where
$M$ is any matrix whose set of columns is~$V$. Clearly, for \emph{rational} matrices $M\in\Q^{m\times n}$ and (finite) sets $V\subseteq\Q^n$ we have $\qsd{M},\qsd{V}\subseteq\Q$. 

The \emph{encoding length} of $\alpha=\tfrac{p}{q}\in \Q$ with $p,q\in\Z$ relatively prime  is
\begin{equation*}
	\enc{\alpha}=1+\lceil\log_2(|p|+1)\rceil+\lceil\log_2(|q|+1)\rceil\,.
\end{equation*}
For a rational vector $v\in\Q^n$ and a rational matrix $M\in\Q^{m\times n}$, we define
\begin{equation*}
	\enc{v}=n+\sum_{j=1}^n\enc{v_j}\quad\text{and}\quad\enc{M}=mn+\sum_{i=1}^m\sum_{j=1}^n\enc{M_{i,j}}\,.
\end{equation*}
Moreover, we denote by $\encmax{M}$ the maximum encoding length of any entry in~$M$, as well as by $\encmax{V}$ the maximum encoding length of all components of vectors in the finite set $V\subseteq\Q^n$.

\subsection{Basics}

Most important, every polyhedron~$P\subseteq\R^n$ is \emph{convex}, i.e., for all $x,y\in P$ and $\alpha\in[0,1]$, we have $\alpha x+(1-\alpha)y\in P$ as well. Moreover, polyhedra are topologically closed subsets of~$\R^n$. 

As the solution sets to finite systems of linear inequalities, polyhedra generalize affine subspaces, which are the solution sets to systems of linear equations. The criterion for $Ax=b$ not being solvable via the existence of some multiplier vector $\lambda\in\R^m$ with $\transpose{\lambda}A=\zeroVec{}$ and $\scalProd{\lambda}{b}\ne 0$ generalizes to systems of linear inequalities in the following way (where both parts of the theorem follow easily from each other).

\begin{theorem}[Farkas-Lemma]\label{thm:Farkas}
	For each $A\in\R^{m\times n}$ and~$b\in\R^m$ the following hold:
	\begin{enumerate}
		\item[(i)] Either $Ax\le b$ is solvable or there is some $\lambda\in\RNonNeg^m$ with $\transpose{\lambda}A=\zeroVec{}$ and $\scalProd{\lambda}{b}<0$ (but not both). 
		\item[(ii)] Either $Ax= b,x\ge\zeroVec{}$ is solvable or there is some $\lambda\in\R^m$ with $\transpose{\lambda}A\ge\zeroVec{}$ and $\scalProd{\lambda}{b}<0$ (but not both). 
	\end{enumerate}
\end{theorem}

It turns out that the ``crucial solutions'' to systems of linear inequalities are obtained as the unique solutions of certain regular linear equation systems, whose components are well-known to be expressible  in the following way.
\begin{theorem}[Cramer's rule]\label{thm:Cramer}
	If $A\in\R^{n\times n}$ is regular, then, for every~$b\in\R^n$, the unique solution $x\in\R^n$ with $Ax=b$ is given by
	\begin{equation*}
		x_j=\frac{\det(\col{A}{1},\dots,\col{A}{j-1},b,\col{A}{j+1},\dots,\col{A}{n})}{\det(A)}
		\quad\text{for all }j\in\ints{n}\,.
	\end{equation*}
\end{theorem}

The following  estimates follows from the Leibniz formula for determinants.


\begin{theorem}\label{thm:encQD}
There is a constant $C>0$ such that the estimates 	
	\begin{equation*}
		\enc{\alpha}\le C\cdot n^2\cdot\encmax{M}\quad\text{for all }\alpha\in\qsd{M}
	\end{equation*}
	and
	\begin{equation*}
		\enc{\beta}\le C\cdot n^2\cdot\encmax{V}\quad\text{for all }\beta\in\qsd{V}
	\end{equation*}
	hold for all $M\in\Q^{m\times n}$ and for all finite sets $V\subseteq\Q^n$.
\end{theorem}

\subsection{Polyhedral and Finitely Generated Cones}

A \emph{cone} is a subset $K\subseteq\R^n$ with $\zeroVec{}\in K$ and $\alpha y\in K$ for all $y\in K$ and $\alpha\in\RNonNeg$. A \emph{polyhedral cone} is a polyhedron that is a cone, or, equivalently, a polyhedron $\polyLe{A}{\zeroVec{}}$ for some $A\in\R^{m\times n}$.  

The (\emph{convex}) \emph{conic hull} of a subset $X\subseteq\R^n$ is the cone
\begin{equation*}
	\ccone{X}=\setDef{\sum_{x\in X'}\alpha_x x}{X'\subseteq X, |X'|<\infty, \alpha_x\ge 0\text{ for all }x\in X'}
\end{equation*}
(with $\ccone{\varnothing}=\{\zeroVec{}\}$) of all \emph{conic combinations} of the vectors in~$X$. 
A cone~$K\subseteq\R^n$ is \emph{finitely generated}, if there is a \emph{finite} set~$X\subseteq\R^n$ with $K=\ccone{X}$. Every vector in a conic hull can  be obtained by a conic combination of few generators:
\begin{theorem}[Carath\'{e}odory's Theorem, conic version]\label{thm:caratheodory}
	For each $X\subseteq\R^n$ and $y\in\ccone{X}$ there is a linearly independent subset $X'\subset X$ (in particular: $|X'|\le n$) with $y\in\ccone{X'}$. 
\end{theorem}
The Farkas-Lemma (Part~(ii) of Theorem~\ref{thm:Farkas}) yields a separation theorem for finitely generated cones.
\begin{theorem}\label{thm:sep}
	If $y\not\in\ccone{X}$ for the finite set $X\subseteq\R^n$, then there is some $a\in\R^n$ with 
	\begin{equation*}
		\scalProd{a}{x}\le0<\scalProd{a}{y}\quad\text{for all } x\in\ccone{X}
	\end{equation*}
	 (i.e., $\ccone{X}\subseteq\halbRaum{a}{0}$, but $y\not\in\halbRaum{a}{0}$).
\end{theorem}

The following result implies that every polyhedral cone is finitely generated, which is of utmost importance for the theory of polyhedra. 
\begin{theorem}\label{thm:polyConeFinite}
	For every matrix $A\in\R^{m\times n}$, there is a finite set $X\subseteq(\qsd{A})^n$ with 
	\begin{equation*}
		\polyLe{A}{\zeroVec{}}=\ccone{X}\,.
	\end{equation*}
\end{theorem}
The \emph{polar} of a cone~$K\subseteq\R^n$ is the  convex cone
\begin{equation*}
	\polar{K}=\setDef{a\in\R^n}{\scalProd{a}{x}\le 0\text{ for all }x\in K}\,.
\end{equation*}
The polar of a finitely generated cone $\ccone{X}$ with a finite set $X\subseteq\R^n$ obviously is the polyhedral cone 
\begin{equation}\label{eq:polarFinite}
	\polar{(\ccone{X})}=\setDef{a\in\R^n}{\scalProd{x}{a}\le 0\text{ for all }x\in X}\,.
\end{equation}
From Theorem~\ref{thm:sep} we also obtain
\begin{equation}\label{eq:polarPoly}
	\polar{(\polyLe{A}{\zeroVec{}})}=\cconeOp\{\row{A}{1},\dots,\row{A}{m}\}
\end{equation}
for each $A\in\R^{m\times n}$. 

Moreover, from Theorem~\ref{thm:sep} one  deduces $\ppolar{(\ccone{X})}=\ccone{X}$, from which one finds, by applying Theorem~\ref{thm:polyConeFinite} as well as (two times)~\eqref{eq:polarFinite}, the following reverse statement to Theorem~\ref{thm:polyConeFinite}.

\begin{theorem}\label{thm:FiniteConePoly}
	For every finite set $X\in\R^n$, there is a matrix $A\in(\qsd{X})^{m\times n}$ with 
	\begin{equation*}
		\ccone{X}=\polyLe{A}{\zeroVec{}}\,.
	\end{equation*}
\end{theorem}

\subsection{The Fundamental Structure of Polyhedra}

The \emph{homogenization} of a polyhedron $\polyLe{A}{b}\subseteq\R^n$ (with $A\in\R^{m\times n}$, $b\in\R^m$) is the polyhedral cone
\begin{equation*}
	\homog{\polyLe{A}{b}}=\setDef{(x,\xi)\in\R^n\oplus\R}{Ax-\xi b\le\zeroVec{},\xi\ge 0}\,.
\end{equation*}
Applying Theorem~\ref{thm:polyConeFinite} to~$\homog{\polyLe{A}{b}}$ as well as Theorem~\ref{thm:FiniteConePoly} to the finitely generated cone
\begin{equation*}
	\ccone{\setDef{(x,1)\in\R^n\oplus\R}{x\in X}\cup\setDef{(y,0)\in\R^n\oplus\R}{y\in Y}}
\end{equation*}
for finite sets $X,Y\subseteq\R^n$, one obtains the following representation theorem for polyhedra, where 
\begin{equation*}
	\conv{X}=\setDef{\sum_{x\in X'}\alpha_x x}{X'\subseteq X, |X'|<\infty, \sum_{x\in X'}\alpha_x=1,\alpha_x\ge 0\text{ for all }x\in X'}
\end{equation*}
denotes the \emph{convex hull} of a set~$X\subseteq\R^n$, and 
\begin{equation*}
	S+T=\setDef{s+t}{s\in S, t\in T}
\end{equation*}
is the \emph{Minkowski sum} of $S,T\subseteq\R^n$.

\begin{theorem}[Weyl-Minkowski Theorem]\label{thm:WeylMinkowski}\mbox{ }
	\begin{enumerate}
		\item[(i)] For every $A\in\R^{m\times n}$ and $b\in\R^m$ there are finite sets~$X,Y\subseteq(\qsd{A,b})^n$ with
		\begin{equation}\label{eq:innerRepr}
			\polyLe{A}{b}=\conv{X}+\ccone{Y}\,.
		\end{equation} 
		\item[(ii)] For all finite sets $X,Y\subseteq\R^n$ there are $A\in(\qsd{X\cup Y})^{m\times n}$ and $b\in(\qsd{X\cup Y})^m$ with
		\begin{equation}\label{eq:outerRepr}
			\conv{X}+\ccone{Y}=\polyLe{A}{b}\,.
		\end{equation} 
	\end{enumerate}
\end{theorem}

Thus, every polyhedron can be represented  by   \emph{outer descriptions} (intersection of finitely many affine halfspaces) and by   \emph{inner descriptions} (Minkowski sum of a polytope and a finitely generated cone). Clearly, algebraically an outer description may contain both linear inequalities and linear equations. Depending on the context, one type of description of a polyhedron can be significantly more convenient to deal with than the other. For instance, from inner descriptions one concludes readily that images of polyhedra under linear maps are polyhedra as well (see Section~\ref{subsec:proj}). In turn, the fact that also preimages of polyhedra under linear maps are polyhedra is easy to  prove via outer descriptions. 
From outer descriptions of polyhedra one also finds immediately that  intersections of finitely many polyhedra  are polyhedra. 

A \emph{rational polyhedron} is a polyhedron for which~$A$ and~$b$, or, equivalently, $X$ and~$Y$, in~\eqref{eq:innerRepr} and~\eqref{eq:outerRepr} can be chosen to be rational matrices, vectors, and sets, respectively. 
Denoting, for a rational polyhedron~$P\subseteq\R^n$, by $\encouter{P}$  the  smallest number~$\alpha$ such that there is an outer description $P=\polyLe{A}{b}$ of~$P$ with a rational matrix~$A$ and a rational vector~$b$ with $\encmax{(A,b)}\le\alpha$, and by $\encinner{P}$  the  smallest number~$\beta$ such that there is an inner description $P=\conv{X}+\ccone{Y}$ of~$P$ with a rational finite sets $X,Y\subseteq\Q^n$ with $\encmax{X\cup Y}\le\beta$, we obtain the following result from Theorem~\ref{thm:encQD}. 
\begin{theorem}\label{thm:outInnComplex}
	There is a constant $C>0$ such that
	\begin{equation*}
		\encinner{P}\le C\cdot n^2\cdot\encouter{P}
		\quad\text{and}\quad
		\encouter{P}\le C\cdot n^2\cdot\encinner{P}
	\end{equation*}
	holds for every rational polyhedron $P\subseteq\R^n$.
\end{theorem}
In particular, whenever the rational system $Ax\le b$ (with $A\in\Q^{m\times n}$ and $b\in\Q^m$) has any solution, then it also has a solution whose encoding length is bounded by a polynomial in~$\enc{(A,b)}$.  This shows that the linear programming feasibility problem is contained in the complexity class~$\NP$, and, via Theorem~\ref{thm:Farkas}, also in $\coNP$. (Of course, it is well-known that this problem is even solvable in polynomial time).

 However, the smallest possible number of inequalities $Ax\le b$ in an outer description and the smallest possible cardinalities of the sets~$X$ and~$Y$ in an inner description of a rational polyhedron~$P$ are not bounded polynomially by each other, in general.

The \emph{characteristic cone} (or \emph{recession cone}) of a polyhedron~$P\subseteq\R^n$ is 
\begin{equation*}
	\charCone{P}=\setDef{y\in\R^n}{x+\coneOp\{y\}\subseteq P\text{ for all }x\in P}\,.
\end{equation*}
\begin{theorem}\label{thm:charCone}
	If $\varnothing\ne P=\polyLe{A}{b}=\conv{X}+\ccone{Y}\subseteq\R^n$ is a polyhedron (with $A\in\R^{m\times n}$, $b\in\R^m$, and $X,Y\subseteq\R^n$ finite), then we have
	\begin{equation*}
		\ccone{Y}=\charCone{P}=\polyLe{A}{\zeroVec{}}\,.
	\end{equation*}
\end{theorem}
The \emph{lineality space} of a polyhedron~$P\subseteq\R^n$ is the largest linear subspace 
\begin{equation*}
	\lineal{P}=\charCone{P}\cap(-\charCone{P})
\end{equation*}
contained in~$\charCone{P}$.
\begin{theorem}\label{thm:lineal}
	For a polyhedron $\varnothing\ne P=\polyLe{A}{b}$  (with $A\in\R^{m\times n}$, $b\in\R^m$) we have
	\begin{equation*}
		\lineal{P}=\kernel{A}\,.
	\end{equation*}
\end{theorem}
A non-empty polyhedron~$\varnothing\ne P\subseteq\R^n$ with   lineality space $\lineal{P}=\{\zeroVec{}\}$ is called \emph{pointed}. For most purposes, it is sufficient to consider pointed polyhedra. In fact, if a polyhedron has a non-trivial lineality space~$L\ne\{\zeroVec{}\}$ then, in most contexts, it is sufficient to investigate instead of~$P$ its orthogonal projection to the orthogonal complement of~$L$, which is a pointed polyhedron. 
Therefore, subsequently we will mainly consider pointed polyhedra. For instance,  all polyhedra that are contained in the nonnegative orthant~$\RNonNeg^n$ as well as all polytopes  are  pointed.

\subsection{Faces of Polyhedra}
\label{subsec:faces}

For $a\in\R^n\setminus\{\zeroVec{}\}$ and $\beta\in\R$,  we denote by 
\begin{equation*}
	\hyperEbene{a}{\beta}=\setDef{x\in\R^n}{\scalProd{a}{x}=\beta}
\end{equation*}
 the boundary hyperplane of the affine halfspace $\halbRaum{a}{\beta}$. A (\emph{proper}) \emph{face} of a polyhedron $P\subseteq\R^n$ is the intersection $F=P\cap\hyperEbene{a}{\beta}$ of~$P$ with the boundary hyperplane of some halfspace $\halbRaum{a}{\beta}\supseteq P$ containing~$P$. The face~$F$ is said to \emph{defined} by the inequality $\scalProd{a}{x}\le\beta$ in this case. Additionally, $\varnothing$ and~$P$ itself are considered (\emph{trivial}) faces of~$P$ (defined by $\scalProd{\zeroVec{}}{x}\le -1$ and $\scalProd{\zeroVec{}}{x}\le 0$, respectively). Non-empty faces are particularly important for optimization, because they are the (non-empty) sets of optimal solutions to linear optimization problems over polyhedra. 

The following result (which is a generalization of  Part~(i) of  Theorem~\ref{thm:Farkas}) provides a characterization of those inequalities that are valid for (and thus define faces of) a non-empty polyhedron.
\begin{theorem}\label{thm:affineFarkas}
	Let $A\in\R^{m\times n}$, $b\in\R^m$, $a\in\R^n$, and $\beta\in\R$ with $P=\polyLe{A}{b}\ne\varnothing$. The inequality $\scalProd{a}{x}\le\beta$ is valid for (all $x$ in) $P$ if and only if there is some $\lambda\in\RNonNeg^m$ with $\transpose{\lambda}A=\transpose{a}$ and $\scalProd{\lambda}{b}\le\beta$.
\end{theorem}

Every face~$F\ne\varnothing$ of a polyhedron $P=\polyLe{A}{b}=\conv{X}+\ccone{Y}$ (with $A\in\R^{m\times n}$, $b\in\R^m$, and $X,Y\subseteq\R^n$ finite) is a polyhedron as well with
\begin{equation*}
	F=\conv{X\cap F}+\ccone{Y\cap\charCone{F}}\,.
\end{equation*}
Moreover, with 
\begin{equation*}
		I=\eqSetMit{A}{b}{F}=\setDef{i\in\ints{m}}{F\subseteq\hyperEbene{\row{A}{i}}{b_i}}\,,
\end{equation*}
we have $F=\setDef{x\in P}{\row{A}{I}x=b_I}$, and, conversely, for every $I\subseteq\ints{m}$, the set $F=\setDef{x\in P}{F\subseteq\hyperEbene{\row{A}{i}}{b_i}}$ is a face of~$P$ with $I\subseteq\eqSetMit{A}{b}{F}$.

In particular, every polyhedron has finitely many faces. Partially ordered by inclusion, they form the \emph{face lattice} $\facLat{P}$ of~$P$ (including~$\varnothing$ and~$P$). For each face~$F\ne\varnothing$ of the polyhedron~$P$ we have $\lineal{F}=\lineal{P}$. Hence, the (non-empty) faces of pointed polyhedra are pointed polyhedra as well.

 The maximal elements in~$\facLat{P}\setminus\{P\}$ are called the \emph{facets} of~$P$. A face of~$P$ is a facet if and only if $\dim(F)=\dim(P)-1$ holds (where the \emph{dimension} $\dim(P)$ of a polyhedron is the affine dimension of its affine hull $\aff{P}$).  The minimal elements in~$\facLat{P}\setminus\{\varnothing\}$ for a pointed polyhedron~$P$ are called the \emph{vertices} of~$P$. The vertices of a polyhedron~$P$ are the faces of dimension zero, i.e., the faces that contain exactly one point (which, of course, is also called \emph{vertex}). 

\begin{theorem}\label{thm:verts}
	For a point $v\in P$ in a polyhedron~$P=\polyLe{A}{b}$ (with $A\in\R^{m\times n}$ and $b\in\R^m$) the following statements are pairwise equivalent.
	\begin{enumerate}
		\item The point~$v$ is a vertex of~$P$.
		\item There are no two points $x,x'\in P\setminus\{v\}$ with $v\in\convOp\{x,x'\}$.
		\item\label{thm:verts:eqsys} There is some $I\subseteq\ints{m}$, $|I|=n$ such that $\row{A}{I}$ is regular with $v=\row{A}{I}^{-1}b$.
	\end{enumerate}
\end{theorem}

A pointed polyhedral cone~$K$ has exactly one vertex, namely $\zeroVec{}$, and therefore, one is more interested in the minimal elements of $\facLat{K}\setminus\{\varnothing,\{\zeroVec{}\}\}$, which are called the \emph{extreme rays} of~$K$. The extreme rays of a (pointed) polyhedral cone are its faces of dimension one. The one-dimensional unbounded faces of a general pointed polyhedron are called its \emph{extreme rays} as well, the one dimensional bounded faces are the \emph{edges}. The edges of a pointed polyhedron~$P$ are of the form $\convOp\{v,w\}$ with two vertices $v\ne w$ of~$P$, and every extreme ray of~$P$ can be written as $v+R$ with a vertex~$v$ of~$P$ and an extreme ray~$R$ of~$\charCone{P}$ (which is pointed if~$P$ is pointed). 

An outer description $P=\setDef{x\in\R^n}{A^=x=b^=,A^{\le}x\le b^{\le}}$ of a polyhedron~$P$ is \emph{irredundant} if removing any equation or any inequalitiy from the system results in a different (larger) polyhedron, and turning any inequality in the system into an equation results in a different (smaller) polyhedron. 
\begin{theorem}\label{thm:outerIrred}
	Let~$\varnothing\ne P\subseteq\R^n$ be a polyhedron, $A^{(1)}\in\R^{m^{(1)}\times n}$, $b^{(1)}\in\R^{m^{(1)}}$, $A^{(2)}\in\R^{m^{(2)}\times n}$, and $b^{(2)}\in\R^{m^{(2)}}$ with
	\begin{equation}\label{eq:outerIrred}
		P\subseteq\setDef{x\in\R^n}{A^{(1)}x=b^{(1)},A^{(2)}x\le b^{(2)}}
	\end{equation}
	such that, for all $i\in\ints{m^{(2)}}$, we have $P\not\subseteq\hyperEbene{\row{A^{(2)}}{i}}{b_i}$.
	\begin{enumerate}
		\item[(i)] Equality in~\eqref{eq:outerIrred} holds if and only if 
		\begin{itemize}
			\item $\aff{P}=\setDef{x\in\R^n}{A^{(1)}x=b^{(1)}}$, and
			\item for each facet~$F$ of~$P$, there is an $i\in\ints{m^{(2)}}$ such that $\scalProd{\row{A^{(2)}}{i}}{x}\le b_i$ defines the facet~$F$ of~$P$.
		\end{itemize}
		\item[(ii)] If equality holds in~\eqref{eq:outerIrred}, then the right-hand-side of~\eqref{eq:outerIrred} is an irredundant outer description of~$P$ if and only if 
		\begin{itemize}
			\item $A^{(1)}$ has full row rank, and
			\item the inequalities in~$A^{(2)}x\le b^{(2)}$ define pairwise distinct facets of~$P$.
		\end{itemize}
	\end{enumerate}
\end{theorem} 

An inner description $P=\conv{X}+\ccone{Y}$ of a  polyhedron~$P$ is \emph{irredundant} if removing any 
point from~$X$ or any vector from~$Y$ results in a different (smaller) polyhedron.
\begin{theorem}\label{thm:innerIrred}
	Let~$\varnothing\ne P\subseteq\R^n$ be a pointed polyhedron, and let $X,Y\subseteq\R^n$ be finite sets with
		$X\subseteq P$
and
		$Y\subseteq \charCone{P}$.
	\begin{enumerate}
		\item[(i)] We have $P=\conv{X}+\ccone{Y}$ if and only if 
		\begin{itemize}
			\item $X$ contains all vertices of~$P$, and
			\item $Y$ contains a nonzero vector from each extremal ray of~$\charCone{P}$.
		\end{itemize}
		\item[(ii)] If $P=\conv{X}+\ccone{Y}$ holds, then $\conv{X}+\ccone{Y}$  is an irredundant inner description of~$P$ if and only if 
		\begin{itemize}
			\item $X$ is the set of vertices of~$P$, and
			\item $Y$ contains exactly one nonzero vector from each extremal ray of~$\charCone{P}$.
		\end{itemize}
	\end{enumerate}
\end{theorem} 

As for a pointed polyhedron~$P=\conv{X}+\ccone{Y}$ (with finite sets $X,Y\subseteq\R^n$), $c\in\R^n$, and  $\omega=\max\setDef{\scalProd{c}{x}}{x\in P}$  we have
$\omega=\infty$ if $\scalProd{c}{y}>0$ for some $y\in Y$, and $\omega=\max\setDef{\scalProd{c}{x}}{x\in X}$ otherwise, Theorem~\ref{thm:innerIrred} implies that a bouned linear optimization problem over a pointed polyhedron~$P$ attains its optimum in a vertex of~$P$.

\subsection{Projections of Polyhedra}
\label{subsec:proj}

We mentioned above that the projection~$P=\pi(Q)$ of a polyhedron $Q\subseteq\R^d$ via a linear map $\pi:\R^d\rightarrow\R^{n}$ with $\pi(y)=Ty$ for a matrix $T\in\R^{n\times d}$ is a polyhedron as well, as, for all $X,Y\subseteq\R^d$, we have 
\begin{equation*}
	\pi(\conv{X}+\ccone{Y})=\conv{\pi(X)}+\ccone{\pi(Y)}\,.
\end{equation*}
We say that a face~$F$ of~$Q$ is \emph{$\pi$-compatible} 
if it 
can be defined by an inequality $\scalProd{\transpose{T}a}{y}\le \beta$ (valid for~$Q$) for some $a\in\R^n$.
\begin{theorem}\label{thm:projPoly}
		For a polyhedron $Q\subseteq\R^d$, a linear projection $\pi:\R^d\rightarrow\R^n$, and the polyhedron $P=\pi(Q)$, the map defined via $F\mapsto \pi(F)$ is an isomorphism between the sublattice of~$\mathcal{L}(Q)$ formed by the $\pi$-compatible faces and~$\mathcal{L}(P)$.
\end{theorem}	

An \emph{extended formulation} for a polyhedron~$P\subseteq\R^n$ is an outer description $Q=\polyLe{A}{b}$ of some polyhedron $Q\subseteq\R^d$ along with a linear projection $\pi:\R^d\rightarrow\R^n$ with $\pi(Q)=P$. 
Such an extended formulation of a polyhedron~$P$ can be much simpler (e.g., in terms of the number of inequalities) than any outer description of~$P$ if 
the complexity of the facets of~$P$ is hidden in lower dimensional parts of the face lattice~$\mathcal{L}(Q)$  (see Theorem~\ref{thm:projPoly}). As, due to
\begin{equation*}
	\max\setDef{\scalProd{c}{x}}{x\in P}=\max\setDef{\scalProd{\transpose{T}c}{y}}{y\in Q}\quad\text{(for all $c\in\R^n$)}\,,
\end{equation*}
linear optimization problems over~$P=\pi(Q)$ can be solved by solving linear optimization problems over~$Q$, extended formulations play an important role in modern mathematical optimization (see~\cite{CCZ09a}). 

We conclude the first part of the article by considering the question how to derive an outer description of a polyhedron from an extended formulation. Here, the fundamental result follows via Theorem~\ref{thm:affineFarkas}. 

\begin{theorem}\label{thm:projPolyOuter}
	Let $Q=\polyLe{D}{g}\subseteq\R^d$ be a polyhedron with $D\in\R^{q\times d}$ and $g\in\R^q$,  suppose the linear projection $\pi:\R^d\rightarrow\R^n$ is defined via $\pi(y)=Ty$ for a matrix $T\in\R^{n\times d}$ and all $y\in\R^d$, and let $\overline{T}$ be any matrix whose rows form a basis of~$\kernel{T}$.
	
	 If $L\in\RNonNeg^{m\times d}$ is a matrix whose rows generate the \emph{projection cone}
	\begin{equation*}
		\setDef{\lambda\in\RNonNeg^q}{\transpose{\lambda}(D\transpose{\overline{T}})=\zeroVec{}}
		=
		\cconeOp\{\row{L}{1},\dots,\row{L}{m}\}\,,
	\end{equation*}
	then every $A\in\R^{m\times n}$ with $AT=LD$ satisfies
	\begin{equation*}
		\pi(Q)=\polyLe{A}{b} \cap \pi(\R^d)
	\end{equation*}
	with $b=Lg$.
\end{theorem}	

If, in the situation of Theorem~\ref{thm:projPolyOuter}, the projection $\pi:\R^{d}\rightarrow\R^n$ is the orthogonal projection to the first~$n$ coordinates, then the projection cone is simply 
\begin{equation}\label{eq:projConeStandard}
	\setDef{\lambda\in\RNonNeg^q}{\scalProd{\col{D}{j}}{\lambda}=0\text{ for all }j\in\ints{d}\setminus\ints{n}}\,,
\end{equation}
and~$A$ can be chosen to consist of the first~$n$ columns of~$LD$. If furthermore $n=d-1$ holds, then a finite generating system for the projection cone~\eqref{eq:projConeStandard} can  be obtained from the relation 
\begin{equation*}
	\setDef{\lambda\in\RNonNeg^{r+s}}{\sum_{i=1}^{r+s}\lambda_i\delta_i=0}
	\ =\ 
	\cconeOp\setDef{\unitVec{i}-\unitVec{k}}{i\in\ints{r},k\in\ints{r+s}\setminus\ints{r}}
\end{equation*}
for all numbers $\delta_1,\dots,\delta_r>0>\delta_{r+1},\dots,\delta_{r+s}$ (which one can easily establish by induction on~$r+s$).

\begin{theorem}[Fourier-Motzkin elimination]\label{thm:FM}
	For a polyhedron $Q=\polyLe{D}{g}\subseteq\R^d$ (with $B\in\R^{q\times d}$ and $g\in\R^q$) and the sets
	\begin{align*}
		I^{>0} &=\setDef{i\in\ints{q}}{D_{i,d}>0}\\
		I^{=0} &=\setDef{i\in\ints{q}}{D_{i,d}=0}\\
		I^{<0} &=\setDef{i\in\ints{q}}{D_{i,d}<0}\,,
	\end{align*}
	the polyhedron in~$\R^{d-1}$ arising from~$Q$ by orthogonal projection to the first $d-1$ coordinates is the set of all $x\in\R^{d-1}$ that satisfy the following system:
	\begin{align*}
		\scalProd{\row{D}{i}}{x} &\le g_i &&\text{for all }i\in I^{=0}\\
		\scalProd{D_{k,d}\row{D}{\ell}-D_{\ell,d}\row{D}{k}}{x} &\le D_{k,d}g_{\ell}-D_{\ell,d}g_k &&\text{for all }k\in I^{>0}, \ell\in I^{<0}
	\end{align*}
\end{theorem}

By applying the Fourier-Motzkin elimination method (Theorem~\ref{thm:FM}) iteratedly, one can compute from an outer descriptions of a polyhedron an outer description of its orthogonal projection to any coordinate subspace. Note that the sizes of the descriptions may blow up exponentially (even if, after each iteration, one removes redundant constraints). 

The case of general projections can be solved by Fourier-Motzin elimination as follows, where we adopt the notation from above. 
Let~$\tilde{n}$ be the rank of~$T$, and assume that the submatrix $T_{\ints{\tilde{n}},\ints{\tilde{n}}}$ of~$T$ is regular. Denoting by~$\tilde{T}\in\R^{d\times d}$ the regular matrix that equals~$\row{T}{\tilde{n}}$ on its first~$\tilde{n}$ rows and 
that has the vectors $\unitVec{i}$, for $i=\tilde{n}+1,\dots,d$, in its other rows, and by $\tilde{\pi}:\R^d\rightarrow\R^{\tilde{n}}$ the orthogonal projection to the first~$\tilde{n}$ coordinates, we have $\pi(y)_{\ints{\tilde{n}}}=\tilde{\pi}(\tilde{T}y)$ for all $y\in\R^d$. Thus, with $\tilde{D}=D\tilde{T}^{-1}$ and $\tilde{Q}=\polyLe{\tilde{D}}{g}$, we find 
	$P=\setDef{x\in\pi(\R^d)}{x_{\tilde{n}}\in\tilde{\pi}(\tilde{Q})}$.

One reason for the interest in computing projections of polyhedra is that one can convert inner to outer and outer to inner description by means of such computations. Indeed, it follows from the discussion of homogenization and the concept of polarity that all such conversions can be done by any method that computes, for a given finite set~$X\subseteq\R^n$, a matrix~$A\in\R^{m\times n}$ with $\ccone{X}=\polyLe{A}{\zeroVec{}}$ (see Theorem~\ref{thm:FiniteConePoly}). Denoting by $T\in\R^{n\times d}$ a matrix whose set of columns is~$X$, we have, by definition,
	$\ccone{X}=\pi(Q)$, where $\pi:\R^d\rightarrow\R^n$ is the projection defined via $\pi(y)=Ty$ for all $y\in\R^d$ and $Q=\RNonNeg^d$. Thus, computing an outer description of~$\pi(Q)$ from the outer description $\RNonNeg^d=\polyLe{-\ident{d}}{\zeroVec{}}$ of~$Q$  solves the problem.

\section{Polyhedra and  Integrality}
\label{sec:int}

Integral points in polyhedra play a crucial role in integer linear programming and combinatorial optimization. The fundamental concept here is the \emph{integer hull} 
\begin{equation*}
	\intHull{P}=\conv{P\cap\Z^n}
\end{equation*}
of a polyhedron~$P\subseteq\R^n$. In this second part, we  describe the most fundamental results on integer hulls, which are important, since, on the one hand, they satisfy
\begin{equation*}
	\max\setDef{\scalProd{c}{x}}{x\in P\cap\Z^n}
	=
	\max\setDef{\scalProd{c}{x}}{x\in \intHull{P}}
\end{equation*}
for all $c\in\R^n$, 
while on the other hand, as convex objects, they are much more convenient to deal with than the discrete sets $P\cap\Z^n$. 
In particular we will treat \emph{integral polyhedra}, i.e., polyhedra~$P$ with $\intHull{P}=P$. Identifying a polyhedron~$P$ as integral is particularly pleasant, as in this case integer linear optimization problems over~$P$ can be treated as (continuous) linear optimization problems over~$P$. 

\begin{theorem}\label{thm:ganzZahligPolyOpt}
	For $A\in\Q^{m\times n}$ and $b\in\Q^m$, defining an integral polyhedron~$\polyLe{A}{b}$, one can, for every $c\in\Q^n$, find $x^{\star}\in\polyLe{A}{b}\cap\Z^n$ with
	\begin{equation*}
		\scalProd{c}{x^{\star}}
		=
		\max\setDef{\scalProd{c}{x}}{x\in\polyLe{A}{b}\cap\Z^n}
	\end{equation*}
	(or conclude that no such $x^{\star}$ exists) in time bounded by a polynomial in $\enc{A}+\enc{b}+\enc{c}$.
\end{theorem}

\subsection{Integer Points in  Polyhedra}

Similarly to the parameterization~\eqref{eq:innerRepr} of points in a polyhedron by an inner description, for a \emph{rational} polyhedron~$P$, we can also parameterize the set of all \emph{integral} points in~$P$. 
For a finite set $Y\subseteq\R^n$, we denote by
\begin{equation*}
	\mono{Y}=\setDef{\sum_{y\in Y}\alpha_yy}{\alpha_y\in\N\text{ for all }y\in Y}
\end{equation*}
(with $\N=\{0,1,2,\dots\}$) the submonoid of~ $\R^n$ generated by~$Y$.
\begin{theorem}\label{thm:intPointsPoly}
	For every rational polyhedron~$P\subseteq\R^n$ there are finite sets $X,Y\subseteq\Z^n$ with
	\begin{equation*}
		P\cap\Z^n=X+\mono{Y}
		\quad\text{and}\quad
		\charCone{P}=\ccone{Y}
		\,,
	\end{equation*}
	such that $\encmax{X\cup Y}$ is bounded by a polynomial in~$n$ and $\encinner{P}$ (and thus, by a polynomial in~$n$ and $\encouter{P}$). 
\end{theorem}
A finite set $H\subseteq\Z^n$ is called an \emph{integral Hilbert basis} (\emph{of $\ccone{H}$})  if 
\begin{equation*}
	\ccone{H}\cap\Z^n=\mono{H}
\end{equation*}
holds. The existence statement of the following theorem is a special case of Theorem~\ref{thm:intPointsPoly}.
\begin{theorem}\label{thm:ratPolyConeHilbert}
	Every rational polyhedral cone~$K$ has an integral Hilbert basis; if~$K$ is pointed, then the integral Hilbert basis of~$K$ is uniquely determined. 
\end{theorem}

Moreover, Theorem~\ref{thm:intPointsPoly} 
implies that, whenever the rational system $Ax\le b$ (with $A\in\Q^{m\times n}$ and $b\in\Q^m$) has any integral solution, then it also has an integral solution whose encoding length is bounded by a polynomial in~$\enc{(A,b)}$.  This shows that the integer linear programming feasibility problem is contained in the complexity class~$\NP$ (it is, however, not in~$\coNP$, unless $\NP=\coNP$). 

Finally, from Theorem~\ref{thm:intPointsPoly} one can derive that the integral hulls of \emph{rational} polyhedra are rational polyhedra as well.

\begin{theorem}\label{thm:intHullPoly}
	For each rational polyhedron~$P\subseteq\R^n$, the integer hull~$\intHull{P}$ of~$P$ is a rational polyhedron, for which $\encouter{\intHull{P}}$ is bounded by a polynomial in~$n$ and $\encouter{P}$.
\end{theorem}

The number of facets of~$\intHull{P}$ is, however, in general not bounded polynomially in the number of facets of~$P$. Furthermore, the integer hull of a non-rational polyhedron needs not even be a polyhedron (see, e.g.,  $P=\setDef{(x_1,x_2)\in\R^2}{x_2\le\sqrt{2}x_1,x_1\ge 1}$). 

A generalization of Theorem~\ref{thm:intHullPoly} to \emph{mixed-integer hulls} of rational polyhedra holds as well: If $P\subseteq\R^n$ is a rational polyhedron, then
	$\convOp\setDef{x\in P}{x_J\in\Z^J}$  is a rational polyhedron for every  $J\subseteq\ints{n}$.

We end this section by a crucial criterion for integrality of (pointed rational) polyhedra.
\begin{theorem}\label{thm:integralityPointedPoly}
	A pointed rational polyhedron is integral (i,e,, $\intHull{P}=P$) if and only if all vertices of~$P$ are integral. 
\end{theorem}

\subsection{Total Dual Integrality}
\label{subsec:TDI}

A rational system $Ax\le b$ with $A\in\Q^{m\times n}$ and $b\in\Q^m$ defining a non-empty polyhedron $\polyLe{A}{b}\ne\varnothing$ is called  \emph{totally dual integral} (\emph{TDI}) if, for every $c\in\Z^n$ with $\max\setDef{\scalProd{c}{x}}{Ax\le b}<\infty$, there is some $y^{\star}\in\N^m$ with $\transpose{A}y^{\star}=c$ and
\begin{equation*}
	\scalProd{b}{y^{\star}}=\min\setDef{\scalProd{b}{y}}{\transpose{A}y=c,y\in\RNonNeg^m}\,,
\end{equation*}
i.e., the dual problem to $\max\setDef{\scalProd{c}{x}}{Ax\le b}$ has an integral optimal solution~$y^{\star}$. It is convenient to  consider TDI all rational systems $Ax\le b$ with $\polyLe{A}{b}=\varnothing$ as well.

\begin{theorem}\label{thm:TDIHilbert}
	A system $Ax\le b$ with integral matrix $A\in\Z^{m\times n}$ and rational right-hand-side vector $b\in\Q^m$ is TDI if and only if, for every face~$F\ne\varnothing$ of~$\polyLe{A}{b}$, the set 
	\begin{equation*}
		\setDef{\row{A}{i}}{i\in\eqSetMit{A}{b}{F}}
	\end{equation*}
	is an integral Hilbert basis.
\end{theorem}

Applying Theorem~\ref{thm:ratPolyConeHilbert} to the rational \emph{normal cones} 
\begin{equation*}
	\setDef{c\in\R^n}{\scalProd{c}{x^{\star}}=\max\setDef{\scalProd{c}{x}}{x\in P}\text{ for all }x^{\star}\in F}
\end{equation*}
of the (minimal) faces~$F$ of a rational polyhedron~$P$, one can construct TDI-systems as in the following result. 
\begin{theorem}\label{thm:ratPolyTDI}
	Let $P\subseteq\R^n$ be a rational polyhedron.
	\begin{enumerate}
		\item There is an integral matrix~$A\in\Z^{m\times n}$ and a rational vector~$b\in\Q^m$ with $P=\polyLe{A}{b}$ such that $Ax\le b$ is TDI, where~$b$ can be chosen to be integral if~$P$ is integral.
		\item\label{item:ratPolyTDIIntb} If there is a rational matrix~$A\in\Q^{m\times n}$ and an integral vector~$b\in\Z^m$ with $P=\polyLe{A}{b}$ such that $Ax\le b$ is TDI, then~$P$ is integral.
	\end{enumerate}
\end{theorem}
Thus, every rational polyhedron admits an outer descriptions by a TDI-system with integral (left-hand-side)  coefficient matrix, though, in general, this will not be an irredundant description in the sense of Section~\ref{subsec:faces}. 
Nevertheless, such a description can provide much structural insight, e.g., within the theory of cutting planes for non-integral polyhedra. The algorithmic problem, however,  to decide whether a given rational system $Ax\le b$ is TDI is $\coNP$-complete (even for~$A$ restricted to the set of node-edge incidence matrices, see Section~\ref{subsec:TU}, of undirected graphs). 

Outer descriptions of (integral) polyhedra with integral right-hand-side vectors 
are very important, because due to Theorem~\ref{thm:ratPolyTDI}(\ref{item:ratPolyTDIIntb}) (and the definition of TDI) they
yield strong duality relations for certain integer linear optimization problems. 
\begin{theorem}
	If the system $Ax\le b$ is TDI with  $A\in\Q^{m\times n}$ and an integral right-hand side $b\in\Z^m$, and $P=\polyLe{A}{b}\ne\varnothing$ holds, then, for every $c\in\Z^n$ with $\max\setDef{\scalProd{c}{x}}{x\in P}<\infty$, we have
	\begin{equation}\label{eq:TDIDuality}
		\max\setDef{\scalProd{c}{x}}{Ax\le b,x\in\Z^n}
		=
		\min\setDef{\scalProd{b}{y}}{\transpose{A}y=c,y\ge\zeroVec{},y\in\Z^m}\,.
	\end{equation}
\end{theorem}

\subsection{Total Unimodularity}
\label{subsec:TU}

A matrix $A\in\{-1,0,1\}^{m\times n}$ is \emph{totally unimodular} (\emph{TU}) if all square submatrices of~$A$ have determinant $-1$, $0$, or~$1$. The (version for pointed polyhedra of the) following result  is a consequence of Theorem~\ref{thm:integralityPointedPoly} and Theorem~\ref{thm:verts} (see part~\eqref{thm:verts:eqsys}) obtained by Cramer's rule (Theorem~\ref{thm:Cramer}).  

\begin{theorem}\label{thm:TUIntPoly}
	If $A\in\{-1,0,1\}^{m\times n}$ is totally unimodular and $b\in\Z^m$ is integral, then $\polyLe{A}{b}$ is an integral polyhedron.
\end{theorem}

If~$A\in\{-1,0,1\}^{m\times n}$ is totally unimodular, then so are all submatrices of~$A$,  $(A,\ident{m})$, $(A,-A)$, and $\transpose{A}$, where the latter allows to conclude the following from Theorem~\ref{thm:TUIntPoly}.
\begin{theorem}
	If $A\in\{-1,0,1\}^{m\times n}$ is totally unimodular, then, for every $b\in\Q^m$, the system $A\le b$ is TDI.
\end{theorem}

In particular, if $A\in\{-1,0,1\}^{m\times n}$ is totally unimodular, then the strong duality relation ~\eqref{eq:TDIDuality} holds for all integral $b\in\Z^m$ and $c\in\Z^n$ (for which the respective optimal values are finite).

The following criterion  is extremely useful for establishing total unimodularity of matrices.

\begin{theorem}[Criterion of Ghouila-Houri]\label{thm:GhouilaHouri}
	A matrix $A\in\{-1,0,1\}^{m\times n}$ is totally unimodular if and only if, for every subset  $I\subseteq\ints{m}$ of row indices, there is a partitioning $I=I^+\uplus I^-$ (with $I^+\cap I^-=\varnothing$) such that
	\begin{equation*}
		\sum_{i\in I^+}\row{A}{i}-\sum_{i\in I^-}\row{A}{i}\in\{-1,0,1\}^n
	\end{equation*}
	holds. (Clearly, a similar characterization via all subsets of column indices holds.)
\end{theorem}

From Theorem~\ref{thm:GhouilaHouri} one readily deduces that every matrix with entries from $\{-1,0,1\}$ that has at most one positive and at most one negative entry per  column is totally unimodular. In particular, the \emph{node-arc incidence matrix} $\inzDir{D}\in\{-1,0,1\}^{V\times A}$ of a  directed graph $D=(V,A)$ (with $A\subseteq V\times V$), defined via
\begin{equation*}
	\inzDir{D}_{v,a}=
	\begin{cases}
		-1 & \text{if }a=(v,w) \text{ for some }w\in V\\
		1 & \text{if }a=(u,v) \text{ for some }u\in V\\
		 0 & \text{otherwise}
	\end{cases}
\end{equation*}
for all $v\in V$ and $a\in A$, is totally unimodular. Thus, the set
\begin{equation*}
	\setDef{x\in\R^A}{\inzDir{D}x=\zeroVec{},\ell\le x\le u}
\end{equation*}
of \emph{circulations} in~$D$ respecting the \emph{integral} lower and upper bounds $\ell,u\in\Z^A$ is an integral polytope. 

The \emph{node-edge incidence matrix} $\inz{G}\in\{0,1\}^{V\times E}$ of an undirected graph $G=(V,E)$ is definied  via
\begin{equation*}
	\inz{G}_{v,e}=
	\begin{cases}
		1 & \text{if }e=\{v,w\}\text{ for some }w\in V\\
		 0 & \text{otherwise}
	\end{cases}
\end{equation*}
for all $v\in V$ and $e\in E$. Theorem~\ref{thm:GhouilaHouri} yields that the node-edge incidence matrix of a graph~$G=(V,E)$ is totally unimodular if and only if the graph is bipartite, i.e., there is a partitioning $V=S\uplus T$ (with $S\cap T=\varnothing$) such that $E\subseteq\setDef{\{s,t\}}{s\in S,t\in T}$. In particular, the matching polytope of a bipartite graph $G=(V,E)$ 
(i.e., the convex hull of the characteristic vectors~$\charVec{M}\in\{0,1\}^E$ ---with $\charVec{M}_e=1$ if and only if $e\in M$--- of all \emph{matchings} $M\subseteq E$ ---with $e\cap e'=\varnothing$ for all $e,e'\in M$, $e\ne e'$) 
has
\begin{align*}
	\sum_{e\in E:v\in e}x_e & \le 1 && \text{for all }v\in V\\
	x_e                     & \ge 0 && \text{for all }e\in E
\end{align*}
as an outer description. 

The most important class of totally unimodular matrices is formed by the \emph{network matrices}. 
Such a network matrix~$N\in\{-1,0,1\}^{A_T\times A}$ arises from a directed graph $D=(V,A)$ and a subset $A_T\subseteq A$ of $|V|-1$ arcs that, viewed as undirected edges, form a spanning tree on~$V$,
by setting, for each $a'\in A_T$ and $a=(v,w)\in A$,
\begin{equation*}
	N_{a',a}=
	\begin{cases}
		1 & \text{if the $v$-$w$-path in }A_T\text{ uses }a'\text{ in forward direction}\\
		-1 & \text{if the $v$-$w$-path in }A_T\text{ uses }a'\text{ in backward direction}\\
		0 & \text{otherwise}
	\end{cases}\,.
\end{equation*}
\begin{theorem}\label{thm:networkMatrices}
	Network matrices are totally unimodular. 
\end{theorem}

In fact, network matrices are the crucial building blocks of the whole class of totally unimodular matrices, as every totally unimodular matrix arises from networks matrices and two special totally unimodular $5\times 5$-matrices by certain operations that preserve total unimodularity. From such structural results one can derive polynomial time algorithms for testing matrices~$A$ for total unimodularity (while, as remarked at the end of Section~\ref{subsec:TDI}, testing a system $Ax\le b$ for total dual integrality is $\coNP$-complete). 

\section{Pointers to Literature}
\label{sec:ref}

Whenever possible, we provide, for the results mentioned in the article, pointers to proofs in~\cite{Sch86}, which is also a great source for  historical notes and many more references than listed here. 

The Farkas-Lemma (Theorem~\ref{thm:Farkas}) dates back to work by Farkas as well as by Min\-kow\-ski at the end of the 19th century (see~\cite[Sect.~7.3]{Sch86}). There are several possibilities to prove the theorem. In fact, the Farkas-Lemma is a special case of  more general separation theorems in convex analysis (see, e.g., \cite[Chap.~2]{Rus06}). A particularly nice and purely linear algebraic proof is due to Conforti, Di Summa, and Zambelli~\cite{CDZ07}, who derive the Farkas-Lemma from the corresponding obstruction of the solvability of linear equation systems. 
Theorem~\ref{thm:Cramer} is due to Cramer (1750)  (see~\cite[Sect.~3.1]{Sch86}).
Theorem~\ref{thm:encQD} follows from~\cite[Thm.~3.2]{Sch86} (see also~\cite{Edm67}). Theorem~\ref{thm:caratheodory} is due to  Carath\'{e}odory~\cite{Car1911} from 1911 (see also~\cite[Thm.~7.1]{Sch86}). Theorems~\ref{thm:polyConeFinite}, \ref{thm:FiniteConePoly}, and~\ref{thm:WeylMinkowski} have their origins in the work of Farkas, Minkowski~\cite{Min1896}, and Weyl~\cite{Wey1935} (see also~\cite[Sect.~ 7.2]{Sch86}). The statements on the components of the  vectors and on the entries of the  matrices (as well as Theorem~\ref{thm:outInnComplex}) follow, e.g., from (the proofs in)~\cite[Sect.~10.2]{Sch86}. An elementary proof of Theorem~\ref{thm:polyConeFinite} can be found in~\cite{Kai09} (see also~\cite[Cor.~7.1a]{Sch86}). For proofs of Theorems~\ref{thm:charCone} and~\ref{thm:lineal} (on the characteristic cones and lineality spaces of polyhedra), we refert to~\cite[Sect.~8.2]{Sch86}.  Theorem~\ref{thm:affineFarkas} is Cor.~7.1h in~\cite{Sch86}. For the other statements in Section~\ref{subsec:faces}, see~\cite[Sect.~8.3--8.9]{Sch86}. Theorem~\ref{thm:projPoly} is folklore (we are not aware of any other explicit reference, thus we refer to~\cite{Kai09b}). Theorem~\ref{thm:projPolyOuter} usually is formulated for orthogonal projections to coordinate subspaces only (see, e.g., \cite[Sect.~2.4]{CCZ09}). An explicit proof in the general setting can be found in~\cite{Kai09b}.
The Fourier-Motzkin method (Theorem~\ref{thm:FM}) is treated in~\cite[Sect.~12.2]{Sch86}. The method is due to Motzkin~\cite{Mot36}, where the idea goes back to work of Fourier in the early 19th century. For the algorithmic problem of converting representations of polyhedra, we refer to the survey by Seidel~\cite{Sei04} and to the software system \texttt{polymake} by Gawrilow and Joswig~\cite{GJ00}  (\url{http://www.opt.tu-darmstadt.de/polymake/}).

The proof of Theorem~\ref{thm:ganzZahligPolyOpt} relies on both that  (continuous) linear programs (Kha\-chi\-yan~\cite{Kha79}) as well as systems of linear Diophantine equations  can be solved in polynomial time (see~\cite[Thm.~16.2]{Sch86}). 
For  proofs of Theorems~\ref{thm:intPointsPoly}, \ref{thm:ratPolyConeHilbert}, and~\ref{thm:intHullPoly}, we refer to~\cite[Sect.~16.2--16.4, 17.2]{Sch86}. 
The fact that the integer hull of a rational polyhedron is a rational polyedron (Theorem~\ref{thm:intHullPoly}) is due to Meyer~\cite{Mey74}, 
the notion of \emph{Hilbert bases} has been introduced by Giles and Pulleyblank~\cite{GP79}, 
where the ideas of the proof of Theorem~\ref{thm:ratPolyConeHilbert} date back to Gordan~\cite{Gor1873}.
The concept of \emph{total dual integrality} has been invented by Edmonds and Giles~\cite{EG77,EG84}. See~\cite[Sect.~22.3]{Sch86}
for proofs of Theorems~\ref{thm:TDIHilbert} and~\ref{thm:ratPolyTDI} (the results being due to~\cite{GP79} and Schrijver~\cite{Sch81}). The $\coNP$-hardness of the TDI-property has been established by Ding, Feng, and Zang~\cite{DFZ08}.
Proofs of the results on total unimodularity can be found in~\cite[Chap.~19]{Sch86}. The connection between totaly unimodular matrices and integral polyhedra (Theorem~\ref{thm:TUIntPoly} and a similar characterization of total unimodularity) is due to Hoffman and Kruskal~\cite{HK56}. Theorem~\ref{thm:GhouilaHouri} has been proved by Ghouila-Houri~\cite{Gho62}. The total unimodularity of network matrices (Theorem~\ref{thm:networkMatrices}) is due to Tutte~\cite{Tut65}. 
The decomposition theorem for totally unimodular matrices mentioned after Theorem~\ref{thm:networkMatrices} has been proved by Seymour~\cite{Sey80}. Cunningham and Edmonds~\cite{CE80} derived a polynomial time test for total unimodularity from that theorem, the asymptotically fastest known algorithm is due to Truemper~\cite{Tru90}.


\end{document}